\documentclass{article}
\usepackage{amstex,amstext,amsthm,a4,amssymb}

\newtheorem{theorem}{Theorem}[section]
\newtheorem{lemma}[theorem]{Lemma}
\newtheorem{proposition}[theorem]{Proposition}
\newtheorem{corollary}[theorem]{Corollary}
\newcommand{\Dom}{\operatorname{Dom}}

\newcommand{\supp}{\operatorname{supp}}
\newcommand{\rank}{\operatorname{rank}}

\begin{document}
\title{Holomorphic Morse Inequalities\\ on Covering Manifolds}
\author{Radu Todor, Ionu\c t Chiose}
\date{}

\maketitle
\begin{abstract}
The goal of this paper is to generalize Demailly's asymptotic holomorphic
Morse inequalities to the case of a covering manifold of a compact
manifold. We shall obtain
estimates which
involve Atiyah's ``normalized dimension'' of the square integrable 
harmonic spaces.
The techniques used are those of Shubin who gave a proof for the usual
Morse inequalities in the presence of a group action relying on Witten
ideas. As a consequence we
obtain estimates for the dimension of the square integrable
holomorphic sections of the pull-back of a line bundle on the base
manifold under some mild hypothesis for the curvature.
\end{abstract}

\section{Introduction}

To emphasize the meaning of our result let us consider a compact
projective manifold $X$, $F$ an ample line bundle on $X$, 
${\widetilde X}$ the universal
covering of $X$ and
${\widetilde F}$ the pull-back of $F$ on $\widetilde X$.
Then Koll\'ar's Theorem 6.4 from \cite{Ko} shows that 
$$\dim_{\pi_1(X)}H^{0}_{(2)}(\widetilde X,K_{\widetilde
X}\otimes\widetilde F)=\dim H^0(X,K_X\otimes F)$$ where we denote as
usual $K_N$ the canonical bundle of a manifold $N$.
It folows from Theorem 0.1 in Demailly \cite{Dem} that
$\dim_{\pi_1(X)}H^{0}_{(2)} (\widetilde X,
K_{\widetilde X}\otimes {\widetilde F}^k)$ has polynomial growth of
order the dimension of $X$ as $k\rightarrow\infty$.
We shall generalize this result to the case of a complex analytic
manifold $M$ on which a discrete group
$\Gamma$ acts freely and properly discontinuous such that $X=M/\Gamma$ is
compact and carries a line bundle $F$ 
satisfying Demailly's condition:
\begin{equation}
\int_{X(\le 1)}(ic(F))^n>0
\end{equation}

\noindent
Let us mention that holomorphic Morse inequalities on non-compact
manifolds have been obtained
before on $q$-concave and $q$-convex manifolds (see Marinescu \cite{Ma} and
Bouche \cite{Bou}).

The cohomology groups of $M$ are usualy infinite--dimensional and we cannot
use the usual dimension. The dimension we shall use is the
$\Gamma$--dimension introduced
by Atiyah in \cite{At}. Let us define the $\Gamma$--dimension of certain
subspaces of
$L^2(M,\widetilde E)$, where $E\rightarrow X$ is a hermitian vector
bundle, $\widetilde E\rightarrow M$ is the pull--back of $E$, $X$ (and hence
$M$) is endowed with a Riemannian metric and
$$
L^2(M,\widetilde E)=\{s:M\rightarrow \widetilde E\,\vert\,s\,
\text{is a measurable section},\mit
\int_M\vert s\vert^2dV<\infty\} 
$$
Let $G$ be a closed subspace in $L^2(M,\widetilde E)$ such that 
$L_{\gamma}G\subset G$
where $L_{\gamma}$ is the action of $\Gamma$ on $L^2(M,\widetilde E)$.
Let $U$ be a fundamental
domain for the action of $\Gamma$ and $(\varphi_m)_m$ an orthonormal 
base for the Hilbert space $G$. Then one defines
\begin{equation}
\dim_{\scriptscriptstyle\Gamma}G:=\sum \limits_m \mit
\int_U\vert\varphi_m(x)\vert^2 dV(x)
\end{equation}
It can be shown that the definition of $\dim_{\scriptscriptstyle
\Gamma}G$ does not
depend on the orthonormal base $(\varphi _m)_m$ or on the fundamental
domain $U$ (see Atiyah \cite{At}).

Let $M$ be a complex analytic manifold of complex dimension $n$ on which a
discrete group
$\Gamma$ acts freely and properly discontinuous such that $X=M/\Gamma$
is compact. Since $M/\Gamma$ is compact one can easily see that the
pull-back of a hermitian metric on $X$ is a complete metric
on $M$ which we consider fixed from now on.
Let $E$ be a hermitian holomorphic vector bundle on $X$ and
 $\widetilde E=\pi^\ast E$ its pull-back,
 where $\pi :M\rightarrow X$ is the projection.

\noindent
Let 
$$
\bar\partial_q:C_{0,q}^{\infty}(M,\widetilde E)\rightarrow
C_{0,q+1}^{\infty}
(M,\widetilde E)
$$ 
be the well-known Cauchy--Riemann operator and 
$$
\delta_q:C_{0,q+1}^{\infty}(M,\widetilde E)
\rightarrow C_{0,q}^{\infty}(M,\widetilde E)
$$ 
the formal adjoint of $\bar\partial_q$.
Then $\Delta^{\prime\prime}_q=\bar\partial _{q-1}\delta _{q-1}+\delta_q\bar 
\partial_q$
is an elliptic differential operator.

Let ${\bar{\partial}}_q:L^2_{0,q}(M,\widetilde E)\rightarrow
L^2_{0,q+1}(M,\widetilde E)$
be the weak maximal extension of $\bar\partial_q$
and likewise we denote by the same letter the weak maximal
extensions of 
${\delta}_q$ and ${\Delta}^{\prime\prime}_q$. 
Let us denote by $N^q(\bar{\partial})$ the kernel of 
${{\bar\partial}}_q$, 
by $R^{q-1}(\bar\partial)$ the range of
${{\bar\partial}}_{q-1}$, 
$N^q(\delta)$ the kernel of ${\delta}
_{q-1}:L^2_{0,q}(M,\widetilde E)
\rightarrow L^2_{0,q-1}(M,\widetilde E)$ and by
$N^q(\Delta^{\prime\prime})$ 
the kernel of ${\Delta}^{\prime\prime}_q
:L^2_{0,q}(M,\widetilde E)\rightarrow L^2_{0,q}(M,\widetilde E)$.

\noindent
By basic results of Andreotti and Vesentini \cite{AV} the hilbertian adjoint of
${{\bar\partial}}$ coincides with ${\delta}$,
${\Delta}^{\prime\prime}_q$ is self--adjoint and 
$$
{\cal H}^q_{(2)}(M,\widetilde E):=N^q(\Delta^{\prime\prime})=N^q (\delta)\cap
N^q(\bar{\partial})
$$
where the first equality is the definition of the space of
$L^2$ harmonic forms. 

\noindent
Let $E$ and $F$ be hermitian holomorphic fibre bundles on $X$ of rank $1$ and $r$
respectively, ${\widetilde F}=\pi ^{\ast}F$, ${\widetilde E}=\pi^{\ast}E$.
Let us denote
$D=D'+\bar {\partial}$ the canonical
connection of $E$ and $c(E)=D^2=D'\bar{\partial}+\bar{\partial}D'$
its curvature form.
Also, let 
$$
X(q)=\lbrace x \in X \mid ic(E)\,\text{has}\, q\, \text{negative
 eigenvalues
and}\, n-q\,
\text{positive eigenvalues} \rbrace
$$ 
and
$$
X(\le q)=X(0)\cup X(1)\cup \ldots \cup X(q).
$$
The main theorem of this paper is:

\begin{theorem}
As $k\to\infty$, the following inequalities
hold for every \\$q=0, 1, \dots, n$ :

\noindent
i) The weak Morse inequalities:
$$
\dim_{\scriptscriptstyle\Gamma}{\cal H}_{(2)}^q(M,\widetilde E^k \otimes \widetilde F)\le
r\frac{k^n}{n!}\int _{X(q)}(-1)^q\left(\frac{i}{2\pi}c(E)\right)^n+
o(k^n).$$
ii) The strong Morse inequalities:
$$\sum\limits_{j=0}^{q}(-1)^{q-j}\dim_{\scriptscriptstyle\Gamma}{\cal H}^{j}_{(2)}(M,
\widetilde E^k\otimes
\widetilde F)
\le r \frac{k^n}{n!} \int_{X(\le q)}(-1)^q\left(\frac{i}{2\pi}c(E)
\right)^n+o(k^n).$$
iii) The asymptotic Riemann-Roch formula:
$$\sum\limits_{j=0}^{n}(-1)^{j}\dim_{\scriptscriptstyle\Gamma}{\cal H}^{j}_{(2)}(M,\widetilde
E^k
\otimes\widetilde F)=
r \frac{k^n}{n!} \int_{X}\left(\frac{i}{2\pi}c(E)\right)^n+o(k^n).
$$
\end{theorem}

\noindent
It follows easily

\begin{corollary} 
Let $E$ and $F$ be as above and suppose that $E$ satisfies  
(1).
Then the space of $L^2$ holomorphic sections satisfies
$$\dim_{\scriptscriptstyle\Gamma}H^{0}_{(2)}(M,\widetilde E^k\otimes
\widetilde F)\approx k^n$$
as $k\rightarrow \infty$. In particular the usual dimension of the
space of $L^2$ holomorphic sections of $E^k$ has the same cardinal as
$|\Gamma|$ for large $k$.
\end{corollary}

\noindent
This generalizes the result for the covering of a projective manifold
by T. Napier \cite{Na}.

\noindent
We wish to express our gratitude to Professor V. Iftimie and Dr. G. Marinescu
for their pertinent suggestions and for the support we have received.

\section{$\Gamma$--dimension and Estimates}

Let $M$ be a real Riemann manifold of dimension $n$, $\Gamma$ a discrete group
acting freely and properly discontinuous on $M$ such that
$X=M/\Gamma$ is compact,
$F\rightarrow X$ a hermitian vector bundle of rank $r$ and 
${\widetilde F}\rightarrow M$ is the pull--back of $F$. 
Let $U$ be a fundamental domain for the action of $\Gamma$.
We identify $L^2(M,\widetilde F)\cong L^2\Gamma \otimes
L^2(U,\widetilde F)\cong L^2\Gamma\otimes L^2(X,F).$
Let us consider ${\cal A}_{\scriptscriptstyle\Gamma }$ the von Neumann
algebra of bounded  operators on $L^2(M,\widetilde F)$ 
which commute with $\Gamma$. 

\noindent
If $A\in{\cal A}_{\scriptscriptstyle\Gamma}$,
then let $K_A\in {\cal D}^{\prime}(M\times M,\widetilde F\otimes_{M\times M}
\widetilde F)$
its kernel. As $A$ is $\Gamma$--invariant, it follows that $K_A \in
{\cal D}^{\prime}(M\times M/\Gamma ,\widetilde F \otimes_{M\times M}
\widetilde F/\Gamma)$ where the action of $\Gamma$ on $M\times M$ 
is $(x,y)\rightarrow (\gamma x,\gamma y)$. 

\noindent
$A\in{\cal A}_{\scriptscriptstyle\Gamma}$ 
is said to be $\Gamma $--Hilbert--Schmidt if $K_A\in L^2(M\times M/\Gamma ,
\widetilde F\otimes_{M\times M}\widetilde F/\Gamma )$ and of
$\Gamma$--trace class if $A=A_1\,A_2$ with 
 $A_1,A_2$ being $\Gamma$--Hilbert--Schmidt.
If $A\in{\cal A}_{\scriptscriptstyle\Gamma}$ is  
of $\Gamma$--trace--class, one can define
\begin{equation}
\operatorname{Tr}_{\scriptscriptstyle\Gamma}A:=\operatorname{Tr}(\varphi A \psi )
\end{equation}
where $\varphi ,\psi\in L^{\infty}_{comp}(M)$ such that $\sum\limits 
_{\gamma\in \Gamma}(\varphi \psi )\circ \gamma =1$. 
If $L\rm\subset L^2(M,\tilde F)$ is
a closed, $\Gamma$--invariant subspace, that is $L$ is a $\Gamma$--module, and
$P_L$ is the ortogonal projection onto $L$, then
\begin{equation}
\dim_{\scriptscriptstyle\Gamma}L:=\operatorname{Tr}_{\scriptscriptstyle\Gamma}
P_L\in[0,\infty ]
\end{equation}

\noindent
This is in short the theory of $\Gamma$--traces. For more results see
Atiyah \cite{At} and Shubin \cite{Sh}. We shall use the following three results;
for the proofs see Shubin \cite{Sh}.

\begin{proposition}\label{prop1} 
Let
$$
0\rightarrow L_0\rightarrow L_1\rightarrow ...\rightarrow L_q
\rightarrow L_{q+1}
\rightarrow\ldots\rightarrow L_n\rightarrow 0  
$$
be a complex of $\Gamma$--modules ($d_q$ commutes with the action
of $\Gamma$ and $d_{q+1}d_q=0$). If
$l_q=\dim_{\scriptscriptstyle\Gamma}L_q
<\infty $ and
$\bar{h}_q=\dim_{\scriptscriptstyle\Gamma}\bar{H}_q(L)$ where 
$$\bar{H}_q(L)=N(d_q)/\overline{R(d_{q-1})}$$ then
\begin{equation}
\begin{split}
\sum\limits_{j=1}^q(-1)^{q-j}\bar{h}_j&\le \sum\limits_{j=1}^q(-1)^{q- 
j}l_j\\
\end{split}
\end{equation}
for every $q=0,1,...,n$ and for $q=n$ the inequality becomes equality.
\end{proposition}

Let $H=H^{\ast}$ be a linear operator in $L^2(M,\widetilde F)$ which commutes  
with the action of $\Gamma$,that is $E_{\lambda}\in{\cal A}
_{\scriptscriptstyle\Gamma}$, where 
$(E_{\lambda})_{\lambda}$ is the spectral family of $H$. Let us denote 
$N_{\scriptscriptstyle\Gamma}(\lambda ,H)=\dim_{\scriptscriptstyle\Gamma}R(E_{\lambda})$ and $h$ the
quadratic form of $H$.

\begin{proposition}\label{prop2}
If $H\ge 0$, then
\begin{multline}
N_{\scriptscriptstyle\Gamma}(\lambda ,H)=\sup\lbrace\dim_{\scriptscriptstyle \Gamma}L\mid
L \;\text{is a}\; \Gamma-\text{module}
\subset\Dom(h),\\ h(f,f)\le \lambda ((f,f)),\forall f\in L\rbrace .
\end{multline}
\end{proposition}

\begin{proposition}\label{prop3}
 If there is $T: L^2(M,\widetilde F)\rightarrow L^2(M,\widetilde F)$
a $\Gamma$--endomorphism (i.e. $T$ commutes with the action of $\Gamma$)  
such that
$((H+T)f,f))\ge\mu ((f,f))$, $f\in\Dom(H)$ and
$\rank_{\scriptscriptstyle\Gamma}T=
\dim_{\scriptscriptstyle\Gamma}\overline{R(T)}\le p$, then
\begin{equation}
N_{\scriptscriptstyle\Gamma}(\mu -\varepsilon ,H)\le p,\;\forall 
\varepsilon >0.
\end{equation}
\end{proposition}

Let $H$ be an elliptic differential operator, formally self--adjoint
of order $2m$ on $\widetilde F$, which commutes with the action of
$\Gamma $.
We shall denote by the same letter $H$ the weak maximal extension of
$H$. If $H$ is strongly elliptic, then
$H$ is bounded from below (see Shubin\cite{Sh}).

\begin{theorem}\label{thm1}
 Let $H_0$ be the self--adjoint operator in 
$L^2(U,\widetilde F\mid_U)$ defined by the restriction of $H$ to $U$ 
with Dirichlet boundary conditions.
($H_0$ is bounded from below and has compact resolvent). Then
\begin{equation}
N_{\scriptscriptstyle\Gamma}(\lambda,H)\ge N(\lambda,H_0),\;\forall
\lambda\in\mathbb R
\end{equation}
where $N(\lambda ,H)=\dim R(F_\lambda)$ if $(F_{\lambda})$ is the 
spectral family of $H_0$.
\end{theorem}

\begin{proof} Let $(e_i)_i$ be an orthonormal basis of
$L^2(U,\widetilde F)$ which consists of eigenfunctions of $H_0$ 
corresponding to the eigenvalues $(\lambda _i)_i$; if we let 
$\widetilde e_i=0$ on $M\setminus U$ and $\widetilde  
e_i=e_i$ on $U$, then $\widetilde e_i\in\Dom(h)$ and
$(L_{\gamma}\widetilde e_i)_{i,\gamma}$
is an orthonormal basis of $L^2(M,\widetilde F)$ and $\widetilde
e_{i,\gamma}=L_{\gamma}\widetilde  e_i\in 
\Dom(h)$. We have $h(\widetilde e_{i,\gamma},\widetilde
e_{i^{'},\gamma^{'}})=\delta_{i,i^{'}}\delta_{\gamma,\gamma^{'}}\lambda_i$.
Let $\Phi_{\lambda}^{0}$
be the subspace spanned by $(\widetilde e_i)_{\lambda_i\le \lambda}$ and
$\Phi_{\lambda}$ the closed subspace spanned by $ (\widetilde e_{i,\gamma} 
)_{
\lambda_i\le \lambda}$. Then 
$$
\dim_{\scriptscriptstyle\Gamma}\Phi_{\lambda}=\sum
((P_{\Phi_{\lambda}}\widetilde e_i,\widetilde
e_i))=\sum\limits_{\lambda_i\le\lambda}
((\widetilde e_i,\widetilde
e_i))=\dim\Phi_{\lambda}^0=N(\lambda,H_0).
$$ 
If $f$ is a linear
combination of $\widetilde e_{i,\gamma},\lambda_i \le \lambda$, then
$h(f,f)
\le\lambda\| f\|^2$ and, as $\Dom(h)$ is complete, we obtain that 
$\Phi_{\lambda}\subset\Dom(h)$ and $h(f,f)\le\lambda\|f\|^2$, 
$f\in \Phi_{\lambda}$. From Proposition \ref{prop2}  it follows that 
$N_{\scriptscriptstyle\Gamma}(\lambda ,H)\ge N(\lambda ,H_0)$.
\end{proof}

Let $s>0$, $U_s=\{x\in M\mid d(x,U)<s\}$ where $d$ is the distance on 
$M$ associated to the Riemann metric on $M$ and $U_{s,\gamma}:=\gamma
U_s$. 
Let $\varphi^{(s)}\in C^{\infty}_0(M)$, $\varphi^{(s)}\ge 0$, 
$\varphi^{(s)}=1$ on $\bar U$ and $\supp\varphi^{(s)}\subset U_s$,
$\varphi^{(s)}_{\gamma}=\varphi^{(s)}\circ \gamma^{-1}$. 
Put 
$$
C^{(s)}_{\gamma}=\frac{\varphi_{\gamma}^{(s)}}
 {\left(\sum\limits_{\gamma}(\varphi^{(s)}_{\gamma})^2\right)^{\frac{1}{2}}}\in
 C^{\infty}_0(M)
$$ 
so that $\sum\limits_{\gamma\in\Gamma}(C^{(s)}_{\gamma})^2=1$.
If $m=1$ (that is $H$ is of order $2$) then
\begin{equation}\label{simb}
H= \sum\limits_{\gamma \in \Gamma}C_{\gamma}^{(s)}HC^{(s)}_{\gamma}-
\sum\limits_{\gamma \in \Gamma}\sigma_0(H)(dC_{\gamma}^{(s)})
\end{equation}
where $\sigma_0$ is the principal symbol of $H$ (see Shubin\cite{Sh}).

Let us assume that $X$ is a complex analytic manifold, $\pi :M\rightarrow X$,
$E$ and $F$ hermitian holomorphic vector bundles on $X$, 
${\widetilde  E}=\pi^{\ast}E$, $\widetilde F=\pi^{\ast}F$. 
Let $\Delta^{\prime\prime}_{k,q}$ be the Laplace--Beltrami
operator on $\Lambda^{0,q}T^{\ast}M\otimes\widetilde E^k\otimes\widetilde F$. 
If $s=k^{-\frac{1}{4}}$,
$H=\frac{1}{k}\Delta^{\prime\prime}_{k,q}$ in \eqref{simb}
then it follows that there is a constant $C$ such that
\begin{equation}
 \frac{1}{k}\Delta^{\prime\prime}_{k,q}\ge \sum\limits_{\gamma\in\Gamma}\frac{1}{k}J 
^{(k)}_{\gamma}
 \Delta^{\prime\prime}_{k,q}J^{(k)}_{\gamma}-\frac{C}{\sqrt{k}}
\operatorname{Id}
 \end{equation}
where $J^{(k)}_{\gamma}=C^{(k^{-\frac{1}{4}})}_{\gamma}$, 
$k\in{\mathbb N}^\ast$.
We have used that $\sigma_0(\Delta^{\prime\prime})(dJ)=\mid
\bar{\partial}J\mid^2\,\operatorname{Id}$  
if $J\in C^{\infty}(M,\mathbb R)$.

\noindent
Let us denote by $V=\Lambda^{0,q}T^{\ast}M\otimes\widetilde E^k\otimes
 \widetilde F$, $H=\frac1k\Delta^{\prime\prime}_{k,q}$ and 
$H^{(k)}_0=\frac{1}{k}\Delta^{\prime\prime}_{k,q}
\mid U_{k^{-\frac{1}{4}}}$ the operator defined in
$L^2(U_{k^{-\frac{1}{4}}},V)$ by the restriction of $\frac{1} 
{k}\Delta^{\prime\prime}_{k,q}$
to $U_{k^{-\frac{1}{4}}}$ with Dirichlet boundary conditions.
Let $(E_{\lambda}^{(k)})_{\lambda}$ be the spectral family of $ H^{(k)}_0$.
In the sequel we fix $\lambda$ and consider $M^{(k)}$ a real number such 
that 
$M^{(k)}\ge \lambda-\inf\operatorname{spec}(H^{(k)}_0)$, where 
$\operatorname{spec} (H^{(k)}_0)$ is the spectrum of $H^{(k)}_0$ to
the effect that 
$$H^{(k)}_0+M^{(k)}E^{(k)}_{\lambda}\ge \lambda\operatorname{Id}.$$ 

\noindent
Define 
\begin{align*}
G^{(k)}_{\gamma}&:L^2(M,V)\rightarrow L^2(M,V)\\
G^{(k)}_{\gamma}&=J^{(k)}_{\gamma}L_{\gamma}M^{(k)}E^{(k)}_{\lambda}
L^{-1}_{\gamma}J_{\gamma}^{(k)},\\
\intertext{(i.e. we trunk the section over $U_{s,\gamma}$, transport it on 
$U_s$, apply the spectral projection and then send it back to 
$U_{s,\gamma}$) and}
G^{(k)}&=\sum\limits_{\gamma\in\Gamma}G^{(k)}_{\gamma}.
\end{align*}
We have
\begin{equation}\label{ineg}
\begin{split}
H+G^{(k)}&\ge\sum\limits_{\gamma\in\Gamma}\left(
J^{(k)}_{\gamma}HJ^{(k)}_{\gamma}
+J^{(k)}_{\gamma}L_{\gamma}M^{(k)}E^{(k)}_{\lambda}L_{\gamma}^{-1}
J_{\gamma}^{(k)}\right)-\frac{C}{\sqrt{k}}\operatorname{Id}\\
&=\sum\limits_{\gamma\in\Gamma}J_{\gamma}^{(k)}L_{\gamma}(H^{(k)}_0+ 
M^{(k)}
E^{(k)}_{\lambda})L_{\gamma}^{-1}J_{\gamma}^{(k)}-\frac{C}{\sqrt{k}}
\operatorname{Id}\\
&\ge\sum\limits_{\gamma\in\Gamma}J_{\gamma}^{(k)}L_{\gamma}\lambda 
L_{\gamma}^{-1}
J_{\gamma}^{(k)}-\frac{C}{\sqrt{k}}\operatorname{Id}\\
&=\left(\lambda -\frac{C}{\sqrt{k}}\right)\operatorname{Id}.
\end{split}
\end{equation}

\begin{lemma}\label{lem1}
$\rank_{\scriptscriptstyle\Gamma}G^{(k)}\le N(\lambda,H^{(k)}_0)$.
\end{lemma}
\begin{proof}The operator 
\begin{gather*}
L^{(1)}_{\alpha}:\bigoplus_{\gamma\in\Gamma} L^2(U_{s,\gamma},V)
\longrightarrow \bigoplus_{\gamma\in\Gamma}L^2(U_{s,\gamma},V),\\
L^{(1)}_{\alpha}((w_{\gamma})_{\gamma}) =
(w_{\alpha^{-1}\gamma})_{\gamma}
\end{gather*}
is a unitary operator for any $\alpha\in\Gamma$. 
Consider $i: L^2(M,V)\rightarrow \bigoplus_{\gamma\in\Gamma}
L^2(U_{s,\gamma},V)$, $i(u)=(u\mid_{U_{s,\gamma}})_{\gamma}$. Then 
$\| u\| \le \| i(u)\|\le C_1 \| u \|$ and hence
$i$ is into and bounded. Moreover $L^{(1)}_{\alpha}i=i\,L_{\alpha}$,
for $\alpha\in\Gamma $.
Let 
\begin{gather*}
F: \bigoplus_{\gamma\in\Gamma} L^2(U_{s,\gamma},V)\rightarrow L^2(M,V),\\
F((w_{\gamma})_{\gamma}) =\sum\limits_{\gamma\in\Gamma}w_{\gamma}.
\end{gather*} 
$F$ is onto, bounded
and $ F\,L^{(1)}_{\alpha}=L_{\alpha}F$, $\alpha\in\Gamma$. 
We define
\begin{gather*}
\widetilde G^{(k)}:\bigoplus_{\gamma\in\Gamma} L^2(U_{s,\gamma},V)\rightarrow
\bigoplus_{\gamma\in\Gamma}L^2(U_{s,\gamma},V),\\
\widetilde G^{(k)}((w_{\gamma})_{\gamma}) =
(J_{\gamma}^{(k)}L_{\gamma}^{-1}M^{(k)}E_{\lambda}^{(k)}L_{\gamma}J
_{\gamma}^{(k)}w_{\gamma})_{\gamma}.
\end{gather*}
Then $\widetilde G^{(k)}$ is bounded, 
commutes with $L^{(1)}_{\alpha}$ and
$G^{(k)}=F\,\widetilde G^{(k)}\,i.$
We define also the operator
\begin{gather*}
K: \bigoplus_{\gamma\in\Gamma}L^2(U_{s,e},V)\rightarrow
\bigoplus_{\gamma\in\Gamma}
L^2(U_{s,\gamma},V),\\
K((w_{\gamma})_{\gamma})=(L_{\gamma}w_{\gamma})_{\gamma}.
\end{gather*}
which is unitary and $K\,L^{(2)}_{\alpha}=L^{(1)}_{\alpha}K$, 
$\alpha\in\Gamma$
where 
\begin{gather*}
L^{(2)}_{\alpha}:\bigoplus_{\gamma\in\Gamma}L^2(U_{s,e},V)\rightarrow
\bigoplus_{\gamma\in\Gamma}L^2(U_{s,e},V),\\
L^{(2)}_{\alpha}((w_{\gamma})_{\gamma})=
(w_{\alpha^{-1}\gamma})_{\gamma}.
\end{gather*} 
Finally, let 
\begin{gather*}
\bar G^{(k)}:
\bigoplus_{\gamma\in\Gamma}L^2(U_{s,e},V)\rightarrow 
\bigoplus_{\gamma\in\Gamma}L^2(U_{s,e},V),\\
\bar G^{(k)}((w_{\gamma})_{\gamma})=
(J_{e}^{(k)}M^{(k)}E_{\lambda}^{(k)}J_{e}^{(k)}w_{\gamma})_{\gamma}.
\end{gather*}
Then 
\begin{equation}\label{comut}
K\,\bar G^{(k)}=\widetilde G^{(k)}K.
\end{equation}
 As $G^{(k)}=F\,
\widetilde G^{(k)}\,i$,
we have that $\rank_{\scriptscriptstyle\Gamma}G^{(k)}\le\rank_
{\scriptscriptstyle\Gamma}\widetilde G^{(k)}$.
The operator $K$ being unitary  
it follows from \eqref{comut}  that 
$\rank_{\scriptscriptstyle\Gamma}\widetilde
G^{(k)}=\rank_{\scriptscriptstyle
\Gamma}\bar G^{(k)}$.
But $R(\bar G^{(k)})$ is closed
because $R(\bar G^{(k)}_{e})$ is closed ($ \bar G^{(k)}_{e}$ is the 
component of $\bar G^{(k)}$ on $L^2(U_{s,e},V)$ and
has finite rank).
If we identify $\bigoplus_{\gamma\in\Gamma}L^2(U_{s,e},V)$ with
$L^2\Gamma\otimes L^2(U_{s,e},V)$ and consider $\bar G^{(k)}$ as an operator
in $L^2\Gamma\otimes L^2(U_{s,e},V)$, then $R(\bar G^{(k)})$ corresponds
to $L^2\Gamma \otimes R(\bar G^{(k)}_{e})$ and hence 
$$
\rank_
{\scriptscriptstyle\Gamma}\bar G^{(k)}=
\rank\bar G^{(k)}_{e}\le\rank E^{(k)}_{\lambda}=N(\lambda,H^{(k)}_{0}).
$$
Now the conclusion follows from the inequality
$$\rank_{\scriptscriptstyle\Gamma}G^{(k)}\le\rank_{\scriptscriptstyle
\Gamma}\widetilde G^{(k)}=\rank 
_{\scriptscriptstyle\Gamma}\bar G^{(k)}.$$
\end{proof}

\begin{proposition}\label{prop4}  
There is a constant $C\ge 0$ such that
\begin{equation}
N_{\Gamma}\left(\lambda ,\frac1k\Delta^{\prime\prime}_{k,q}\right)\le N\left(\lambda 
 +
\frac{C}{\sqrt{k}},\frac1k\Delta^{\prime\prime}_{k,q}\mid
U_{k^{-\frac{1}{4}}}\right)
\;\lambda\in{\mathbb R},\; k\in{\mathbb N}^\ast
\end{equation}
\end{proposition} 

\begin{proof} Proposition \ref{prop3} with $\mu =\lambda -\frac{C}{\sqrt 
{k}}$ and $p=N\left(\lambda,\frac1k\Delta^{\prime\prime}
_{k,q}\mid U_{k^{-\frac{1}{4}}}\right)$,
\eqref{ineg} and Lemma \ref{lem1} entail 
$$
N_{\scriptscriptstyle\Gamma}
\left(\lambda-\frac{C}{\sqrt{k}}-\varepsilon ,\frac1k\Delta^{\prime\prime}_{k,q}\right)
\le N\left(\lambda ,\frac1k\Delta^{\prime\prime}_{k,q}  \mid U_{\frac{1}{\sqrt[4]{k} 
}}\right),\,\varepsilon >0
$$
Replacing $\lambda$ with $\lambda+\frac{C}{\sqrt{k}}+\varepsilon$,
we obtain
$$ N_{\scriptscriptstyle\Gamma}\left(\lambda
,\frac{1}{k}\Delta^{\prime\prime}_{k,q}
\right)\le
N\left(\lambda +\frac{C}{\sqrt{k}}+\varepsilon ,\frac{1}{k}\Delta ^{\prime\prime}_{k,q}
\mid U_{k^{-\frac{1}{4}}}\right)$$
When $\varepsilon \rightarrow 0$ it follows
$$ 
N_{\scriptscriptstyle\Gamma}\left(\lambda,\frac1k\Delta^{\prime\prime}_{k,q}\right)\le
N\left(\lambda +\frac{C}{\sqrt{k}},\frac1k\Delta^{\prime\prime}_{k,q}\mid U
_{k^{-\frac{1}{4}}}\right)
$$
\end{proof}

\section{Holomorphic Morse Inequalities}

Let $M$ be a Riemannian manifold of dimension $n$ with volume element 
$d\sigma$.
Let $E$ and $F$ be hermitian vector bundles on $M$, $\rank E=1$,
$\rank F=r$, with $D$ and
$\nabla$ the canonical connections, $S$ a continuous section in
$\Lambda^1_{\mathbb R}T^{\ast}M\otimes_{\mathbb
R}\operatorname{Hom}_{\mathbb C}(F,F)$ and $V$
a continuous section in $\operatorname{Herm}(F)$. Let
$\nabla_k$ be the connection in $E^k\otimes F$. We denote the endomorphisms
$\operatorname{Id}_{E^k}\otimes S$ and $\operatorname{Id}_{E^k}\otimes V$ by 
$S$ and $V$. Given $\Omega\Subset M$, let
\begin{gather}
Q_{\Omega ,k}(u)=\int_{\Omega}\left(\frac{1}{k}\mid\nabla_ku+Su
\mid^2-(Vu,u)\right)
d\sigma\,,\\ 
\Dom(Q_{\Omega,k})=W^1_0(\Omega,E^k\otimes F)\notag
\end{gather}
where by $W^1_0$ we denote the Sobolev space.
Let $V_1(x)\le ... \le V_r(x)$ be the eigenvalues of V(x). We shall
use the following

\begin{theorem}[Demailly \cite{Dem}] 
The counting function of the eigenvalues of 
$Q_{\Omega,k}$ satisfies for every $\lambda\in\mathbb R$ 
the following asymptotic estimates as $k\rightarrow \infty$:
\begin{multline}
\sum\limits_{j=1}^r\int_{\Omega}\nu_B(V_j+\lambda)d\sigma \le
\liminf\;k^{-\frac{n}{2}}N(\lambda,Q_{\Omega,k})\le\limsup\;
k^{-\frac{n}{2}}N(\lambda,Q_{\Omega,k})\le \\
\sum\limits_{j=1}^r\int_{\overline\Omega}\bar{\nu}_{B}(V_j+\lambda)
d\sigma
\end{multline}  
where $B$ is the magnetic field of the connection $D$ and
\begin{equation}
\nu_B(\lambda)=\frac{2^{s-n}\pi^{-\frac{n}{2}}}{\Gamma\left(
\frac{n}{2}-s+1\right)}\,B_1\dotsm B_s\sum\limits_{(p_1,\dotsc,p_n)\in
{\mathbb N}^s}\left[
\lambda-\sum\limits_{j=1}^s(2p_j+1)B_j\right]_+^{\frac{n}{2}-s}
\end{equation}
if $B_1(x)\ge ...\ge B_s(x)$ are the absolute values of the non--zero 
eigenvalues of $B$, $[\lambda]^0_+=0$ for $\lambda\le 0$, $[\lambda]^ 
0_+=1$ for $\lambda >0$ and $\bar{\nu}_B=\lim\;\nu_B(\lambda+\varepsilon ), 
\varepsilon\searrow 0$.
\end{theorem}


We assume now that $M$ is a hermitian complex analytic manifold of complex 
dimension $n$, $E$ and $F$ hermitian holomorphic vector bundles,
$F_k=E^k\otimes F$ and $\Delta^{\prime\prime}_k$ the Laplace--Beltrami operator
on $F_k$. If $\alpha_1(x),\dotsc,\alpha_n(x)$ are the eigenvalues of $ic(E)(x)$
with respect to the metric on $M$ then, as in Demailly \cite{Dem}, we deduce that there is
a countable set $A\subset\mathbb R$ such that
\begin{equation}
\lim\,k^{-n}N\left(\lambda,\frac{1}{k}\Delta^{\prime\prime}_{k,q}\mid\Omega\right)=
r\sum\limits_{\mid J\mid =q}\int_{\Omega}\nu_B(2\lambda+\alpha_{C(J)}-\alpha
_{J})d\sigma
\end{equation}
for $\lambda\in{\mathbb R}\setminus A$,
where $\alpha_J=\sum\limits_{j\in J}\alpha_j$, $C(J)=\{1,...,n\} 
\setminus J$.

Let $M$ be an analytic complex manifold of dimension $n$ and $\Gamma$
a discrete group which acts freely and properly discontinuous on M
such that $X=M/\Gamma$ is compact. We choose a hermitian metric on X 
and we lift it on $M$. Let $E$ and $F$
be hermitian holomorphic vector bundles on $X$, $\rank E=1$, $\rank F=r$ and
$\widetilde E=\pi^{\ast}E$, $\widetilde F=\pi^{\ast}F$.

\noindent
Let $E(\cdot,\Delta^{\prime\prime}_{k,q})$ the spectral family of
the self-adjoint
operator $\Delta^{\prime\prime}_{k,q}$ in $L^2_{0,q}
(M,\widetilde E^k\otimes\tilde F)$
and $L^{\lambda,k}_{q}=R\left(E\left([0,\lambda ],\frac{1}{k}\Delta 
^{\prime\prime}_{k,q}\right)\right)$.
Then 
$$
E([0,\lambda
],\Delta^{\prime\prime}_{k,q}){\bar{\partial}}_{q-1}=
{\bar{\partial}}_{q-1}E([0,\lambda ],\Delta^{\prime\prime}
_{k,q-1})
$$ 
on $L^{\lambda,k}_{q-1}$ and it follows that $\bar{\partial}_{q-1}L^ 
{\lambda,k}_{q-1}\subset
L^{\lambda,k}_{q}$. 
If $\bar{\partial}_{q}^{\lambda}$ denotes the restriction
$$\bar{\partial}_{q}:L^{\lambda,k}_{q}\rightarrow L^{\lambda,k}_{q+1}$$ 
 then
\begin{equation}\label{harm}
\{u\in L^{\lambda,k}_{q}\mid \bar{\partial}^{\lambda}_{q}u=0,
(\bar{\partial}^{\lambda}_{q-1})^{\ast}u=0\}=
{\cal H}^{q}_{(2)}(M,\widetilde E^k\otimes \widetilde F)
\end{equation}
(see Shubin \cite{Sh}).
By definition
$N_{\scriptscriptstyle\Gamma}(\lambda,\frac{1}{k}\Delta^{\prime\prime}_{k,q})
=\dim_{\scriptscriptstyle\Gamma}
L^{\lambda,k}_{q}$.
There is a complex
\begin{equation}
0\rightarrow L^{\lambda,k}_{0}\rightarrow L^{\lambda,k}_{1}\rightarrow\dotsm
\rightarrow L^{\lambda,k}_{n}\rightarrow 0
\end{equation}
From Proposition \ref{prop1} 
we get
$$ \sum\limits_{j=1}^{q}(-1)^{q-j}\dim_{\scriptscriptstyle\Gamma}
(N(\bar{\partial}^{\lambda}_{q})
/\overline{R(\bar{\partial}^{\lambda}_{q-1})})\le \sum\limits_{j=1}^{q}(-1)^{ 
q}
N_{\scriptscriptstyle\Gamma}\left(\lambda,\frac{1}{k}\Delta^{\prime\prime}
_{k,q}\right)$$
for $q=0,1,\dotsc,n$ and for $q=n$ the inequality becomes equality.

From \eqref{harm}, 
$N(\bar{\partial}^{\lambda}_{q})
/\overline{R(\bar{\partial}^{\lambda}_{q-1})}
\simeq {\cal H} ^{q}_{(2)}(M,\widetilde E^k\otimes \widetilde F)$. 
From Theorem \ref{thm1}
and Proposition \ref{prop4} it follows that there is a constant $C$ such that
\begin{equation}
N(\lambda,\frac{1}{k}\Delta^{\prime\prime}_{k,q}\mid U)\le N_{\scriptscriptstyle\Gamma}(\lambda,
\frac{1}{k}\Delta ^{\prime\prime}_{k,q})\le
N\left(\lambda+\frac{C}{\sqrt{k}},\frac{1}{k}\Delta^{\prime\prime}_{k,q}\mid U
_{\frac{1}
{\sqrt[4]{k}}}\right)
\end{equation}
But
$$N(\lambda,\frac{1}{k}\Delta^{\prime\prime}_{k,q}\mid
U)=r\,k^n\sum\limits_{\mid J\mid=q}
\int_U\nu_B(2\lambda+\alpha_{C(J)}-\alpha_{J})d\sigma +o({k^n})$$
and
\begin{multline*}
\limsup\; k^{-n}N\left(\lambda+\frac{C}{\sqrt{k}},\frac{1}{k}\Delta
^{\prime\prime}_{k,q}
\mid U_{k^{-\frac{1}{4}}}\right)\le
\lim\,k^{-n}N(\lambda+\varepsilon 
,\frac{1}{k}
\Delta^{\prime\prime}_{k,q}\mid U_{\varepsilon})\\
=r\sum\limits_{\mid J\mid =q}\int_{U_{\varepsilon}}\nu_B(2\lambda+ 
2\varepsilon+
\alpha_{C(J)}-\alpha_{J)})d\sigma
\end{multline*}
so when $\varepsilon \rightarrow 0$ we get
$$
\lim\,k^{-n}N\left(\lambda+\frac{C}{\sqrt{k}},\frac{1}{k}\Delta^{\prime\prime}
_{k,q}\mid
U_{k^{-\frac{1}{4}}}\right)=r\sum\limits_{\mid J\mid=q}\int_{\bar U}
\bar{\nu}_B(2\lambda+\alpha_{C(J)}-\alpha_{J})d\sigma
$$
for every $\lambda\in{\mathbb R}\setminus A$. As $\partial U=\bar U\setminus U$
is of measure zero because $U$ is a fundamental domain, it follows that
$$
N_{\scriptscriptstyle\Gamma}(\lambda,\frac{1}{k}\Delta^{\prime\prime}_{k,q})
=r\,k^n\sum\limits_{\mid J\mid=q}\int_{U}\nu_B(2\lambda+\alpha_{C(J)}-
\alpha_{J})d\sigma +o(k^n)
$$
for $\lambda\in{\mathbb R}\setminus A$.
Hence for $\lambda \rightarrow 0$, $\lambda\in{\mathbb R}\setminus A$ we obtain
\begin{equation}\label{weak}
\sum\limits_{j=0}^{q}(-1)^{q-j}\dim_{\scriptscriptstyle\Gamma}
{\cal H}^j_{(2)}(M,\widetilde E^k\otimes
\widetilde F)\le k^n\sum\limits_{j=0}^{q}(-1)^{q-j}I^j+o(k^n)
\end{equation}
where
$$I^j=r\sum\limits_{\mid J\mid=q}\int_{U}\nu
_B(\alpha_{C(J)}-\alpha_{J})d\sigma$$
which is calculated in Demailly \cite{Dem} :
\begin{equation}\label{int}
I^j=\frac{r}{n!}\int_{M(j)\cap
U}(-1)^j\left(\frac{i}{2\pi}c(\widetilde E)
\right)^n=\frac{r}{n!}\int_{X(j)}(-1)^j\left(\frac{i}{2\pi}c(E)\right)^n
\end{equation}
with $M(j)=\{x\in M\mid\;ic(\widetilde
E)(x)\,\text{has}\,j\,\text{negative
eigenvalues and}
\,n-j\,\text{positive ones}\}$. 
We have used that $c(\widetilde E)$ is the lifting of $c(E)$.
Theorem 1.1 now follows from \eqref{weak} and \eqref{int}.

\bibliographystyle{plain}

\vskip2cm
\noindent
Address of the authors:  Department of Mathematics,\\
University of Bucharest,
Str. Academiei 14,
70109 Bucharest, Romania.

\noindent
e-mail: radone@@skylab.math.unibuc.ro, chiose@@skylab.math.unibuc.ro

\end{document}